%% file: main.tex
\documentclass[11pt,reqno]{amsart}
\usepackage[utf8]{inputenc}
\usepackage[T1]{fontenc}
\usepackage{lmodern}
\usepackage{amsfonts,amsthm,amsmath,amssymb}
\usepackage{graphicx}
\usepackage[dvipsnames,table]{xcolor}   
\usepackage{enumerate}
\usepackage{hyperref}
\usepackage{tikz-cd}
\usepackage[margin=1in]{geometry}
\usepackage[
    maxbibnames=99,
    backend=biber,
    style=alphabetic,
    sorting=nyt,
    doi=false,
    giveninits=true
]{biblatex}
\DeclareFieldFormat{pages}{#1}
\renewbibmacro{in:}{%
  \ifentrytype{article}
    {}
    {\bibstring{in}%
     \printunit{\intitlepunct}}}
\DeclareFieldFormat
[article,inbook,incollection,inproceedings,patent,thesis,unpublished]
  {title}{\mkbibemph{#1}}
\DeclareFieldFormat{journaltitle}{#1\isdot}
\DeclareFieldFormat[article]{volume}{\mkbibbold{#1}}
\DeclareFieldFormat[article]{number}{\bibstring{number}\addnbspace #1}

\renewbibmacro*{journal+issuetitle}{%
  \usebibmacro{journal}%
  \setunit*{\addspace}%
  \iffieldundef{series}
    {}
    {\newunit
     \printfield{series}%
     \setunit{\addspace}}%
  \printfield{volume}%
  \setunit{\addspace}%
  \usebibmacro{issue+date}%
  \setunit{\addcomma\space}%
  \printfield{number}%
  \setunit{\addcolon\space}%
  \usebibmacro{issue}%
  \setunit{\addcomma\space}%
  \printfield{eid}
  \newunit}

\usetikzlibrary{arrows.meta}

\newtheoremstyle{mystyle}
  {}
  {}
  {\itshape}
  {}
  {\bfseries}
  {.}
  { }
  {\thmname{#1}\thmnumber{ #2}\thmnote{ (#3)}}

\theoremstyle{mystyle}
\newtheorem{Thm}{Theorem}

\newtheorem{Cor}[Thm]{Corollary}
\newtheorem{Prop}[Thm]{Proposition}

\theoremstyle{definition}

\newtheorem{Ex}[Thm]{Example}

\theoremstyle{remark}
\newtheorem{Rmk}[Thm]{Remark}

\newcommand{\Z}{\mathbb{Z}}
\newcommand{\Q}{\mathbb{Q}}

\newcommand{\Sk}{\mathcal S}
\newcommand{\Int}{\mathrm{int}}

\author{Qiuyu Ren}
\address{Department of Mathematics, University of California, Berkeley, Berkeley, CA 94720, USA}
\email{\href{mailto:qiuyu_ren@berkeley.edu}{qiuyu\_ren@berkeley.edu}}

\author{Fang-Rong Zhan}
\address{Department of Mathematics, North Carolina State University, Raleigh, NC 27695, USA}
\email{\href{mailto:fzhan@ncsu.edu}{fzhan@ncsu.edu}}

\title{Insensitivity of Khovanov homology under rim surgery}

\begin{document}

\begin{abstract}
We show that rim surgery on a smooth, oriented, properly embedded surface in $B^4$ does not change the map on Khovanov homology induced by the surface, answering a question raised by Hayden and Sundberg \cite{hayden2024khovanov}. More generally, the same insensitivity holds for a class of local surgery operations that we call local annulus replacements. Combined with the work of Hayden--Sundberg, this result yields exotic pairs of surfaces in $B^4$ that cannot be related by any sequence of local annulus replacements. We give two proofs, both making essential use of Khovanov skein lasagna modules.
\end{abstract}

\maketitle

\section{Introduction}
Rim surgery, introduced by Fintushel and Stern \cite{fintushel1997surfaces}, is a procedure for producing new embedded surfaces in $4$-manifolds from a given one. 
Let $\Sigma$ be a smooth, oriented, properly embedded surface in a compact oriented smooth $4$-manifold, and let $\gamma\subset\Sigma$ be a simple closed curve along which a normal framing of $\Sigma$ is chosen. Around $\gamma$, we can use the normal framing along $\gamma$ to identify the manifold-surface pair with $(S^1\times B^3,S^1\times B^1)$. We then replace each untangled arc $B^1$ with a (1,1)-tangle $T_K\subset B^3$ obtained by cutting open a knot $K\subset S^3$, i.e., we remove $(S^1\times B^3,S^1\times B^1)$ and reglue $(S^1\times B^3,S^1\times T_K)$. The process is called \textit{rim surgery} along $\gamma$ with pattern $K$.

Several sufficient conditions are known under which rim surgery
preserves the topological isotopy type. For example, this holds when the surface complement is simply connected \cite{fintushel1997surfaces}. It also holds when $\gamma$ bounds a locally flat topological disk in the surface complement and the rim surgery is performed using the normal framing induced by that disk with one additional twist \cite{juhasz2021transverse}. Thus, rim surgery provides a way to construct potentially exotic pairs of surfaces. To detect changes in smooth isotopy type, previous studies have used Seiberg--Witten invariants \cite{fintushel1997surfaces}, Ozsv\'ath--Szab\'o mixed invariants \cite{mark2013knotted}, or perturbed link Floer homology \cite{juhasz2021transverse}.

In \cite{hayden2024khovanov}, Hayden and Sundberg produce, for every integer $g\ge0$, examples of exotic pairs of genus-$g$ surfaces in $B^4$ bounding a common knot in $S^3$ that are distinguishable by their induced maps on Khovanov homology \cite{khovanov2000categorification}. However, it was unknown whether the Khovanov cobordism map could distinguish exotic surfaces related by rim surgery, as pointed out in \cite{hayden2024khovanov}. In the present work, we answer this question negatively and prove a more general result.

\begin{Thm}\label{thm:main}
Let $\Sigma$ and $\Sigma'$ be smooth, oriented, properly embedded surfaces in $B^4$ with common boundary $L\subset S^3$. If $\Sigma$ and $\Sigma'$ are related by a local annulus replacement, then, for
every commutative ring $R$, the maps
\[
R\cong Kh(\emptyset;R)\to Kh(L;R)
\]
induced by $\Sigma$ and $\Sigma'$ agree.
\end{Thm}

Here, by a \textit{local annulus replacement} we mean the operation of replacing a local pair $(S^1\times B^3,S^1\times B^1)$ around a curve $\gamma\subset\Sigma$ by $(S^1\times B^3,A)$ where $A$ is any smoothly embedded annulus with boundary $S^1\times\partial B^1$. Thus, rim surgery is a special case of a local annulus replacement.

When $\Sigma$ and $\Sigma'$ are disks, a local annulus replacement operation reduces to a connected sum operation with a local $2$-knot, which is well known to be undetectable by the Khovanov cobordism map; thus, Theorem~\ref{thm:main} is only interesting for surfaces with positive genus. (Conversely, a connected sum with a local $2$-knot can be realized by a local annulus replacement that happens near a point on the surface.)

For comparison, the analogous insensitivity of the link Floer homology $\widehat{HFL}$ under rim surgeries was established as a special case of \cite{juhasz2021transverse}.

Applying Theorem~\ref{thm:main} to Hayden--Sundberg's examples of exotic surfaces \cite[Theorem~1]{hayden2024khovanov} gives the following corollary.
\begin{Cor}\label{cor:main}
For every integer $g\ge0$, there exists an exotic pair $\Sigma,\Sigma'$ of smooth, oriented, properly embedded genus-$g$ surfaces in $B^4$ with common boundary a knot in $S^3$ such that $\Sigma$ and $\Sigma'$ are not related by any sequence of local annulus replacements. In particular, they are not related by any sequence of rim surgeries.
\end{Cor}

In Sections~\ref{sec:proof1} and~\ref{sec:proof2}, we present two different proofs of Theorem~\ref{thm:main}. In both proofs, the first observation is that both maps in Theorem~\ref{thm:main} factor through the Khovanov skein lasagna module of $S^1\times B^3$ rel the boundary link $S^1\times\partial B^1$. The difference is that the first proof\footnote{The first proof only applies when the coefficient ring $R$ is a field of characteristic not equal to $2$.} uses the Khovanov skein lasagna module with $0$-dimensional inputs, as originally introduced by Morrison--Walker--Wedrich \cite{morrison2022invariants}, while the second proof, significantly simpler, uses the Khovanov skein lasagna module with $1$-dimensional inputs, introduced by Ren--Sullivan--Wedrich--Willis--Zhang \cite{ren2025khovanov}. In Section~\ref{sec:generalization}, we mention two natural generalizations of Theorem~\ref{thm:main}.

\section*{Acknowledgments}
We thank Tye Lidman, Robert Lipshitz, Ciprian Manolescu, Gheehyun Nahm, and Ian Zemke for helpful discussions, and Tye Lidman for comments on an earlier draft of this paper. The second author would also like to thank Jiahuang Chen, Cai Zeng, Daren Chen, and Mihai Marian for sharing their thoughts. The research was partially carried out during the 2026 IAS/PCMI Graduate Summer School and Research Program. We are grateful to the organizers for fostering a stimulating environment for productive discussion and research.

This research was partially conducted during the period QR served as a Clay Research Fellow. FZ was partially supported by the NSF grant DMS-2506277.

\section{First proof}\label{sec:proof1}
In this proof, we need to assume the coefficient ring $R$ is a field $k$ of characteristic not equal to $2$. We assume the reader has some familiarity with Khovanov skein lasagna modules \cite{morrison2022invariants}; see \cite[Section~2]{ren2024khovanov} for an exposition.

We will use $\Sk_0^2(X;L)=\Sk_0^2(X;L;k)$ to denote the Khovanov skein lasagna module of the pair $(X,L)$ with $k$-coefficients, where $X$ is an oriented compact smooth $4$-manifold and $L\subset\partial X$ is a framed oriented link. By convention, $\Sk_0^2$ takes the Khovanov--Rozansky $\mathfrak{gl}_2$ homology $KhR_2(\bullet):=KhR_2(\bullet;k)$ as the input link homology theory. The theory $KhR_2$ is defined for framed oriented links and framed oriented link cobordisms. Since $KhR_2$ is functorially isomorphic to $Kh$ up to mirroring the link and applying grading shifts determined by the framing, we may fix normal framings on $L,\Sigma,\Sigma'$ so that $\Sigma,\Sigma'$ are related to each other by a local framed annulus replacement along some curve $\gamma\subset\Sigma$, and prove Theorem~\ref{thm:main} with $Kh$ replaced by $KhR_2$.

For later convenience, we may assume that the framing on $\Sigma$ is chosen so that the normal framing along $\gamma$ given by the normal of $\gamma$ in $\Sigma$ followed by the normal framing of $\Sigma$ along $\gamma$ is the even framing. If $\gamma$ is null-homologous on $\Sigma$, this is automatic; otherwise, we may change the framing on $\Sigma$, if necessary, to achieve this.

Let $\nu(\gamma)$ be a closed tubular neighborhood of $\gamma$ and pick an identification $(\nu(\gamma),\nu(\gamma)\cap\Sigma)\cong(S^1\times B^3,S^1\times B^1)$ compatible with the normal framing of $\Sigma$. Write $A=\Sigma'\cap(S^1\times B^3)$. As a shorthand, let $1_2$ denote the $2$-component link $S^1\times\partial B^1$ in $S^1\times S^2$ equipped with the orientation and framing induced from $S^1\times B^1=\Sigma\cap(S^1\times B^3)$. Viewed in $S^1\times S^2$, $1_2$ is null-homologous and carries the standard framing.\smallskip

\textbf{Step 1}: Reduce to the special case of rim surgeries.

Both maps $KhR_2(\emptyset)\to KhR_2(L)$ induced by $\Sigma$ and $\Sigma'$ factor as
\begin{equation}\label{eq:factors_0}
\Sk_0^2(\emptyset;\emptyset)\to\Sk_0^2(S^1\times B^3;1_2)\to\Sk_0^2(B^4;L),
\end{equation}
where the second arrows in the two factorizations agree for $\Sigma$ and $\Sigma'$ because the surfaces agree outside $S^1\times B^3=\nu(\gamma)$. The domain and the target of \eqref{eq:factors_0} are canonically identified with $KhR_2(\emptyset)$ and $KhR_2(L)$ respectively, and Theorem~\ref{thm:main} amounts to showing that the skeins $\Sigma$ and $\Sigma'$ represent the same element in $\Sk_0^2(B^4;L)$. We would thus like to understand the two elements in $\Sk_0^2(S^1\times B^3;1_2)$ represented by the skeins $S^1\times B^1$ and $A$, as well as the map in \eqref{eq:factors_0} out of $\Sk_0^2(S^1\times B^3;1_2)$. To this end, we need the following crucial computation from Manolescu--Walker--Wedrich \cite{manolescu2023skein}, which is only available over a field.

\begin{Prop}[{\cite[Section~4.4]{manolescu2023skein}}]
If $k$ is a field, we have $\Sk_0^2(S^1\times B^3;1_2)\cong k^4$ concentrated in (homological, quantum) bidegrees $(0,0),(0,0),(0,2),(1,4)$. Moreover, one can fix an identification $\Sk_{0,0,0}^2(S^1\times B^3;1_2)\cong k^2$, such that for any framed oriented knot $K$, the skein $S^1\times T_K$ represents the element $(1,n(K))\in k^2\cong\Sk_{0,0,0}^2(S^1\times B^3;1_2)$ for some integer $n(K)$. Here, $T_K$ is the $(1,1)$-tangle associated to $K$.
\end{Prop}

The integer $n(K)$ (really, an integral multiple of $1\in k$) is determined as follows. The proof of Proposition~4.15 in \cite{manolescu2023skein} shows that the Khovanov--Rozansky $\mathfrak{gl}_2$ chain complex of $K$ is chain homotopy equivalent to a direct sum of shifted copies of one-term complexes $k[X]/X^2$ and two-term complexes $k[X]/X^2\xrightarrow{X}k[X]/X^2$. Counting the number of $k[X]/X^2\xrightarrow{X}k[X]/X^2$ summands with signs (depending on whether the first $k[X]/X^2$ appears in odd or even homological grading), one obtains the integer $n(K)$. (Analogously, it is easy to check, using the fact that $J_K(1)=2$, that the signed count of $k[X]/X^2$ summands is exactly $1$, which gives the first coordinate in $(1,n(K))$.)

\begin{Ex}\label{ex:trefoil}
For the crossingless diagram of the $0$-framed unknot $U$, $CKhR_2(U)=k[X]/X^2$ where the unit $1\in k[X]/X^2$ has degree $(0,-1)$. Thus $n(U)=0$, and $[S^1\times B^1]=(1,0)\in\Sk_{0,0,0}^2(S^1\times B^3;1_2)$.

If $K=T(2,3)$ is the $0$-framed right-handed trefoil knot, then $KhR_2(K)\cong k^4$ is concentrated in degrees $(0,1),(0,3),(-2,5),(-3,9)$ (this requires $k$ to have characteristic not equal to $2$). By grading considerations, we may conclude that $CKhR_2(K)$ is chain homotopy equivalent to the direct sum of one shifted copy of $k[X]/X^2$ where $1\in k[X]/X^2$ is in degree $(0,1)$, and one shifted copy of $k[X]/X^2\xrightarrow{X}k[X]/X^2$ where $1\in k[X]/X^2$ in the first term is in degree $(-3,7)$. Consequently, upon fixing a sign convention, $n(T(2,3))=1$, and we may write $[S^1\times T_{T(2,3)}]=(1,1)\in\Sk_{0,0,0}^2(S^1\times B^3;1_2)$.
\end{Ex}

We continue with the proof of Theorem~\ref{thm:main}. We claim that, in fact, any annular skein $A$ in $S^1\times B^3$ rel $1_2$ represents an element of the form $(1,m(A))\in\Sk_{0,0,0}^2(S^1\times B^3;1_2)$ for some $m(A)\in k$. This can be seen by a cap-off argument. Capping $(S^1\times B^3;1_2)$ off in $S^4$ by gluing in $B^2\times S^2$ with the disjoint union of two disks bounding $1_2$, one dotted and one undotted, yields a map
\begin{equation}\label{eq:capoff}
k^2\cong\Sk_{0,0,0}^2(S^1\times B^3;1_2)\to\Sk_{0,0,0}^2(S^4)\cong k.
\end{equation}
The image of $[A]$ under this gluing map is the evaluation of a once-dotted sphere using the functor $KhR_2$. Every once-dotted sphere in $S^4$ evaluates to $1$ under $KhR_2$, and applying this fact to the two annuli $S^1\times B^1$ and $S^1\times T_{T(2,3)}$ in Example~\ref{ex:trefoil} shows that the linear map \eqref{eq:capoff} is given by $(a,b)\mapsto a$. Then, we may conclude that any annulus $A$ in $S^1\times B^3$ rel $1_2$ represents a class in $\Sk_{0,0,0}^2(S^1\times B^3;1_2)\cong k^2$ of the form $(1,m(A))$, as desired.

Consequently, by linearity, it suffices to prove Theorem~\ref{thm:main} when $\Sigma'$ is obtained from $\Sigma$ by a rim surgery (or indeed by a rim surgery with pattern the trefoil knot) using the framing of $\Sigma$. We therefore assume that $A=S^1\times T_K$ for some framed oriented knot $K$.\smallskip

\textbf{Step 2}: The case of a rim surgery.

We make a topological reduction, which may not be necessary for the proof but simplifies the picture. Since $\gamma$ is contractible in $B^4$, we may find a smoothly embedded oriented disk $D\subset\Int(B^4)$ bounding $\gamma$. By general position, we may assume that the interior of $D$ intersects $\Sigma$ transversely; in particular, $D\cap\Sigma$ is the disjoint union of $\gamma$ with finitely many oriented points in $\Int(D)$.

We claim that the disk $D$ may be chosen so that along $\gamma$,
\begin{enumerate}[(i)]
\item the normal direction of $\gamma$ in $D$ agrees with the framing of $\Sigma$ along $\gamma$;
\item the normal direction of $\gamma$ in $\Sigma$ extends to a normal framing of $D$.
\end{enumerate}
The obstructions to both (i) and (ii) take values in $\Z$. Performing a positive boundary twist on $D$ along $\gamma$ increases the obstruction to (i) by $1$, while changing the obstruction to (ii) $\pm1$, with the sign depending on the side of $\gamma$ on $\Sigma$ on which the new intersection point between $D$ and $\Sigma$ appears. A negative boundary twist has the opposite effect. By the even assumption on the normal framing of $\Sigma$ along $\gamma$ at the beginning of the proof, the two obstructions sum to an even integer. Therefore, we may perform some number of boundary twists to make both obstructions vanish, giving (i) and (ii).

Let $B=D_1\times D_2$ denote a small tubular neighborhood of $D$ in $B^4$, where $D_1$ is a slightly expanded copy of $D$ and $D_2$ is a disk normal to $D$. We have arranged that
\begin{itemize}
\item $\Sigma$ and $\Sigma'$ agree outside the $4$-ball $B$;
\item $\partial B\cap\Sigma=\partial B\cap\Sigma'=H_{2,n}$, where $H_{2,n}$ denotes the parallel cable of the $0$-framed Hopf link $H=\partial D_1\times\{*\}\cup\{*\}\times\partial D_2$ in which the first component of $H$ is $2$-cabled, with oppositely oriented components, and the second component of $H$ is $n$-cabled, with various orientations.
\item $B\cap\Sigma$ is the disjoint union of $n$ vertical (i.e., parallel to $D_2$) disks bounding the $n$ components in $H_{2,n}$ and a standard $\partial$-parallel annulus bounding the two components in $H_{2,n}$. $B\cap\Sigma'$ is the disjoint union of $n$ vertical disks bounding the $n$ components in $H_{2,n}$, together with the spun annulus of $T_K$, which bounds the two components of $H_{2,n}$. This spun annulus is contained in a neighborhood of $\partial D_1\times D_2$ in $B$, identified with $S^1\times B^3$, where it appears as a copy of $S^1\times T_K$.
\end{itemize}
It remains to show that the surfaces $B\cap\Sigma$ and $B\cap\Sigma'$ in $B$ induce identical maps $KhR_2(\emptyset)\to KhR_2(H_{2,n})$. We take out a smaller ball $B'$ in $\Int(B)$ that meets $\Sigma$ and $\Sigma'$ each in some shrunk copies of the $n$ vertical disks, and regard the rest of $\Sigma$ and $\Sigma'$ as cobordisms from $H_{0,n}$ to $H_{2,n}$, denoted $A_n(T_U)$ and $A_n(T_K)$, then prove that their induced maps on $KhR_2$ are equal.

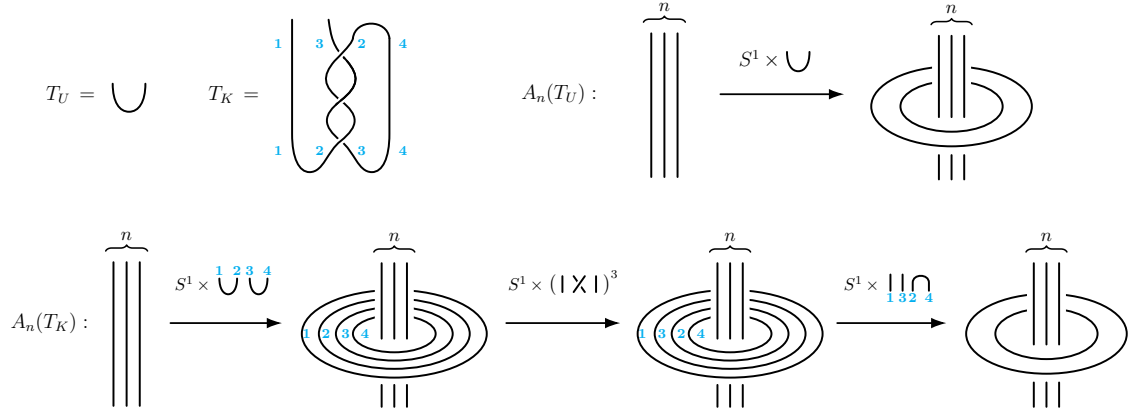
\begin{figure}[htbp]
\centering
\resizebox{1\linewidth}{!}{%
        \input{rim}%
    }
\caption{First row: The tangles $T_U$ and $T_K$ in bridge position, with $T_K$ shown for $K=T(2,3)$, and the cobordism $A_n(T_U)$ given by an annulus birth. Second row: The cobordism $A_n(T_K)$, expressed as the composition of two annulus births, three braidings of the two middle circles, and one annulus death. The $n$ free strands in each frame close up near infinity.}
\label{fig:movie}
\end{figure}

Putting $T_U$ and $T_K$ into bridge position, we obtain the movie descriptions of $A_n(T_U)$ and $A_n(T_K)$ shown in Figure~\ref{fig:movie}. Here, the height function on $B\backslash\Int(B')$ is chosen to be $S^1$-invariant in the $S^1\times B^3$ neighborhood of $\partial D_1\times D_2$ so that its restriction to each $B^3$ slice induces the chosen bridge presentations for $T_U$ and $T_K$. In general, if $T_K$ is in $2b$-bridge position, $A_n(T_K)$ is a composition of $b$ annulus births, some number of braidings among the $2b$ circles, and $b-1$ annulus deaths.

Grigsby--Licata--Wehrli \cite[Section~7]{grigsby2018annular} shows that the braid group action by $B_m$ on the $m$-cable of a framed oriented knot in $S^3$ descends to a symmetric group action by $S_m$ on the Khovanov homology of the cable. More generally, their argument shows that if $K\cup L\subset S^3$ is a framed oriented link where $K$ is a knot, and $K(m)\cup L$ denotes the link obtained from $K\cup L$ by $m$-cabling $K$, then the $B_m$-action on $K(m)\cup L$ permuting the $m$-cable descends to an $S_m$-action on $KhR_2(K(m)\cup L)$. In particular, the case $m=2$ implies that the signs of the braidings in the description of $A_n(T_K)$ do not affect the induced map on $KhR_2$; in other words, we may change crossings in the chosen bridge presentation of $T_K$ without affecting $KhR_2(A_n(T_K))\colon KhR_2(H_{0,n})\to KhR_2(H_{2,n})$. Since any knot $K$ can be changed to the unknot by a sequence of crossing changes, we see that $KhR_2(A_n(T_K))=KhR_2(A_n(T_U))$, concluding the proof of Theorem~\ref{thm:main}.\qed

\section{Second proof}\label{sec:proof2}
Rozansky \cite{rozansky2010categorification} and Willis \cite{willis2021khovanov} defined a theory of Khovanov homology for links in connected sums of $S^1\times S^2$. In \cite{ren2025khovanov}, using Khovanov skein lasagna modules, Ren--Sullivan--Wedrich--Willis--Zhang showed the functoriality of Rozansky--Willis homology with $\Q$-coefficients for framed oriented links in $\#(S^1\times S^2)$ and link cobordisms in a large class of $4$-manifolds called \textit{relative $1$-handlebody complements}, defined as the complement of a (possibly empty or disconnected) $4$-dimensional $1$-handlebody in another $4$-dimensional $1$-handlebody.

More precisely, if $X_0,X_1$ are compact oriented $4$-manifolds with only $0$- and $1$-handles, $X_0\subset\Int(X_1)$, $L_0,L_1$ are links in $\partial X_0,\partial X_1$, then every link cobordism in $X_1\backslash\Int(X_0)$ from $L_0$ to $L_1$ induces a map from the Rozansky--Willis homology of $L_0$ to that of $L_1$, and these maps are functorial under composition of such link cobordisms. The first author subsequently removed the restriction on the coefficients and extended the functoriality to any coefficient ring (and indeed to the chain level up to chain homotopy) \cite[Chapter~4]{ren2026khovanov}.

We clarify the convention. Four versions of Rozansky--Willis homology can be considered: two of which are defined using the bounded-below, unital Rozansky projectors, indicated by a $+$ superscript, and two defined using the bounded-above, counital Rozansky projectors, indicated by a $-$ superscript. It is the bounded-above versions that are functorial with respect to the cobordisms that we described and supply the input link homology theory of Khovanov skein lasagna modules with $1$-dimensional inputs. For the two flavors of bounded-above Rozansky--Willis homology, one is identically zero unless the link is null-homologous, and the other is identically zero unless the link is $\Z/2$-null-homologous and may be nonzero for links that are not integrally null-homologous. For the purpose of this proof, the distinction does not matter, but for concreteness we will work with the former flavor. We follow the convention in \cite[Chapters~4.1.1 and~4.2.1]{ren2026khovanov} (compare also \cite[Section~2.2]{ren2025khovanov}) and use the bounded-above Rozansky--Willis homology for null-homologous links with $R$-coefficients, denoted $KhR_{2,\text{0-div}}^-(\bullet):=KhR_{2,\text{0-div}}^-(\bullet;R)$. When $L$ is a link in $S^3$, we have $KhR_{2,\text{0-div}}^-(L)\cong KhR_2(L)$ canonically. Our second proof of Theorem~\ref{thm:main} essentially uses the functoriality result for Rozansky--Willis homology, although to put it in parallel with the first proof, we write the proof in terms of Khovanov skein lasagna modules with $1$-dimensional inputs, defined in \cite[Section~3.2]{ren2025khovanov} and \cite[Chapter~4.1.1]{ren2026khovanov}. Following the notation in \cite{ren2026khovanov}, we use $\Sk_0^{2,\text{1d},\text{0-div}}(X;L):=\Sk_0^{2,\text{1d},\text{0-div}}(X;L;R)$ to denote the $1$-dimensional-input Khovanov skein lasagna module with input homology theory $KhR_{2,\text{0-div}}^-$.\footnote{In the notation of \cite{ren2025khovanov}, $\Sk_0^{2,\text{1d},\text{0-div}}$ is denoted $\bar\Sk_0^{2,O}$, which appears in Remark~3.5(1).} Like the $0$-dimensional-input Khovanov skein lasagna module, $\Sk_0^{2,\text{1d},\text{0-div}}$ comes with a homological grading and a quantum grading which will be used in the proof.
\begin{Rmk}
\begin{enumerate}
\item Strictly speaking, following \cite{ren2026khovanov}, the functoriality is only proved for the ``mul'' version of Rozansky--Willis homology, and we should have used the notation $KhR_{2,\text{mul},\text{0-div}}^-(L)$ instead of $KhR_{2,\text{0-div}}^-(L)$. The distinction is invisible in the case when the link $L$ is in a single copy of $\#(S^1\times S^2)$ or the empty space (as opposed to a disjoint union of multiple copies of $\#(S^1\times S^2)$). These are the only cases considered here, so we suppress the ``mul'' subscript. See \cite[Chapter~2.6.2.2]{ren2026khovanov} for more discussion.
\item Similarly, the corresponding skein lasagna module should strictly be written $\Sk_0^{2,\text{1d},\text{mul},\text{0-div}}$. We suppress the ``mul'' superscript because the two versions agree when $R$ is a field, whereas for a general commutative ring only the ``mul'' version is defined. Thus, the shortened notation should cause no ambiguity.
\end{enumerate}
\end{Rmk}

After these background explanations, we return to the proof of Theorem~\ref{thm:main}. As in the first proof, consider a small $S^1\times B^3$ neighborhood of the curve along which the local annulus replacement is performed. The two maps $KhR_2(\emptyset)\to KhR_2(L)$ induced by $\Sigma$ and
$\Sigma'$ then factor as
\begin{equation}
\Sk_0^{2,\text{1d},\text{0-div}}(\emptyset;\emptyset)\to\Sk_0^{2,\text{1d},\text{0-div}}(S^1\times B^3;1_2)\to\Sk_0^{2,\text{1d},\text{0-div}}(B^4;L),
\end{equation}
where the second arrows in the factorizations agree for $\Sigma$ and $\Sigma'$. Here, as in the first proof, $1_2$ denotes the union of two oppositely oriented, standardly framed $S^1$ fibers in $S^1\times S^2$. We claim that the first maps are also equal for $\Sigma$ and $\Sigma'$, which would finish the proof of Theorem~\ref{thm:main}.

The key difference here compared to the first proof is that the relevant graded piece $\Sk_{0,0,0}^{2,\text{1d},\text{0-div}}(S^1\times B^3;1_2)\cong KhR_{2,\text{0-div}}^{-,0,0}(1_2)\cong R$ is free of rank $1$ instead of rank $2$. This is a classical computation: The bounded-above Rozansky projector closure on $2n$ strands gives a projective resolution of Khovanov's arc algebra $H^n$ \cite{khovanov2002functor} over its enveloping algebra. Thus, $KhR_{2,\text{0-div}}^-(1_2)$ is isomorphic, as a bigraded $R$-module, to the Hochschild homology of $H^1$ with $R$-coefficients. More precisely, $H^1=R[X]/X^2$, and in homological degree $0$ we see $$\Sk_{0,0,*}^{2,\text{1d},\text{0-div}}(S^1\times B^3;1_2)\cong KhR_{2,\text{0-div}}^{-,0,*}(1_2)\cong HH_0(H^1)=H^1/[H^1,H^1]=H^1=R[X]/X^2$$ is free of rank $2$, concentrated in quantum degrees $0$ and $2$, whence the claim. One can check from the Hochschild description that the standard annulus $S^1\times B^1$ represents $1\in R\cong\Sk_{0,0,0}^{2,\text{1d},\text{0-div}}(S^1\times B^3;1_2)$.

To finish the proof, we apply the cap-off argument again as in the first proof. As before, consider the map 
\begin{equation}\label{eq:capoff_2}
R\cong\Sk_{0,0,0}^{2,\text{1d},\text{0-div}}(S^1\times B^3;1_2)\to\Sk_{0,0,0}^{2,\text{1d},\text{0-div}}(S^4)\cong R
\end{equation}
obtained by gluing in $B^2\times S^2$ with an undotted cap and a dotted cap. The image of any class $[A]$ represented by an annulus $A\subset S^1\times B^3$ is the evaluation of a dotted sphere in $S^4$, which is equal to $1\in R$. By considering $S^1\times B^1$, we conclude that \eqref{eq:capoff_2} is the identity map on $R$. It then follows that every annulus $A$ in $S^1\times B^3$ bounding $1_2$ represents the class $1\in R\cong\Sk_{0,0,0}^{2,\text{1d},\text{0-div}}(S^1\times B^3;1_2)$, proving Theorem~\ref{thm:main}.\qed

\section{Generalizations}\label{sec:generalization}
In this section we mention two natural generalizations of our result on the insensitivity of Khovanov homology under local annulus replacements.

First, if $X$ is any oriented smooth $4$-manifold, one could try to distinguish smooth oriented surfaces in $X$ bounding some $L\subset\partial X$ by studying their class in the skein lasagna module $\Sk_0^2(X;L)$. Our proof of Theorem~\ref{thm:main} applies to show that if $X$ is simply connected and $\Sigma,\Sigma'$ are two such surfaces related by a sequence of local annulus replacements, then they represent the same class in $\Sk_0^2(X;L)$.

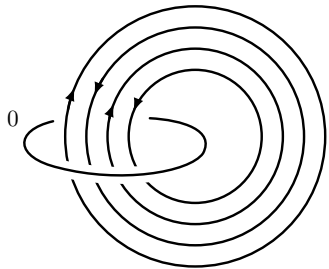
\begin{figure}[t]
\centering\scalebox{.8}{
\begin{tikzpicture}[
    line width=1.15pt,
    line cap=round,
    >={latex[length=3mm,width=2.3mm]}
  ]
\begin{scope}[shift={(1.05,0)}]
  \foreach \r in {1.10,1.45,1.80,2.15}
    \draw (0,0) circle[radius=\r];
  \def\a{158} 
  \def\d{10}  
  \foreach \r/\s in {2.15/1,1.80/-1,1.45/1,1.10/-1}
    \draw[->]
      ({\a+\d*\s}:\r)
      arc[start angle={\a+\d*\s},end angle=\a,radius=\r];
\end{scope}
\def\loop{(-1.30,.25)
  .. controls (-1.95,.10) and (-1.95,-.35) .. (-1.20,-.52)
  .. controls (-.20,-.78) and (.95,-.60) .. (1.18,-.28)
  .. controls (1.38,.10) and (.85,.27) .. (.3,.3)}
\draw[white,line width=8pt] \loop;
\draw \loop;
\node[left] at (-1.72,.30) {$0$};
\end{tikzpicture}}
\caption{A standard diagram of $1_4\subset S^1\times S^2$ with alternating orientations.}
\label{fig:1_2n}
\end{figure}

Second, albeit somewhat impractical, one might try to distinguish two framed smooth oriented surfaces $\Sigma,\Sigma'$ (say in $B^4$) bounding some common framed oriented link $L\subset S^3$ by taking each of their parallel $n$-cables $\Sigma(n),\Sigma'(n)$ bounding the $n$-cable $L(n)$ of $L$, and studying the Khovanov cobordism maps of the cables. Here, all cables are taken with respect to the given normal framings. Our second proof of Theorem~\ref{thm:main} generalizes to show that if $\Sigma$ and $\Sigma'$ are related by a sequence of local annulus replacements, then $Kh(\Sigma(n))=Kh(\Sigma'(n))$ for any $n$. The additional relevant input is that, if $1_{2n}=1_2(n)\subset S^1\times S^2$ denotes the $n$-cable of the link $1_2\subset S^1\times S^2$, then $$KhR_{2,\text{0-div}}^{-,0,*}(1_{2n})\cong HH_0(H^n)\cong H^n/[H^n,H^n],$$ whose quantum degree $0$ part is isomorphic to $\oplus_{a\in B^n}R\cong R^{C_n}$ where $H^n$ denotes Khovanov's arc algebra on $2n$ strands over $R$, $B^n$ is the set of crossingless matchings on $2n$ points, and $C_n=|B^n|=\frac{1}{n+1}\binom{2n}{n}$ is the $n$-th Catalan number. Choose a standard diagram of $1_{2n}$ in which the orientations alternate, as in Figure~\ref{fig:1_2n}. Then the cap-off of $1_{2n}$ in $S^1\times B^3$ given by $(S^1\times B^1)(n)$, the $n$-cable of $S^1\times B^1$, represents the element $1$ living in the $R$-factor in $KhR_{2,\text{0-div}}^{-,0,0}(1_{2n})\cong\oplus_{a\in B^n}R$ that corresponds to the matching of the $2n$ strands given by $(S^1\times B^1)(n)$. The elements $1$ in other $R$-factors are similarly represented by various collections of standard annuli capping off $1_{2n}$. By a cap-off argument using $B^2\times S^2$ with various combinations of dotted and undotted caps, one can show that the $n$-cable of any annulus $A\subset S^1\times B^3$ bounding $1_2$, or indeed any disjoint union of $n$ annuli in $S^1\times B^3$ bounding $1_{2n}$ that pairs the components of $1_{2n}$ in the same way as $(S^1\times B^1)(n)$, represents the same element in $KhR_{2,\text{0-div}}^{-,0,0}(1_{2n})$ as $(S^1\times B^1)(n)$.

\printbibliography

\end{document}

%% file: rim.tex
\usetikzlibrary{arrows.meta,decorations.pathreplacing,calc}

\definecolor{strandblue}{RGB}{0,174,239}

\tikzset{
  every picture/.style={line cap=round,line join=round},
  diagram/.style={draw=black,line width=1.15pt},
  movie arrow/.style={-{Latex[length=3.2mm,width=2.2mm]},line width=1.05pt},
  blue label/.style={text=strandblue,font=\scriptsize\bfseries},
  arc label/.style={text=strandblue,font=\scriptsize\bfseries,inner sep=.35pt},
  pics/cup/.style={
    code={
      \draw[diagram] (-.42,.27)
        .. controls (-.42,-.20) and (-.25,-.45) .. (0,-.45)
        .. controls (.25,-.45) and (.42,-.20) .. (.42,.27);
    }
  },
  pics/doublecup/.style={
    code={
      \draw[diagram] (-.95,.22)
        .. controls (-.95,-.18) and (-.79,-.40) .. (-.60,-.40)
        .. controls (-.41,-.40) and (-.25,-.18) .. (-.25,.22);
      \draw[diagram] (.25,.22)
        .. controls (.25,-.18) and (.41,-.40) .. (.60,-.40)
        .. controls (.79,-.40) and (.95,-.18) .. (.95,.22);
      \node[blue label] at (-.95,.52) {1};
      \node[blue label] at (-.25,.52) {2};
      \node[blue label] at (.25,.52) {3};
      \node[blue label] at (.95,.52) {4};
    }
  },
  pics/crossblock/.style={
    code={
      \draw[diagram] (-.72,-.34)--(-.72,.34);
      \draw[diagram] (.72,-.34)--(.72,.34);
      \draw[diagram] (.24,-.34)--(-.24,.34);
      \draw[diagram,preaction={draw=white,line width=3.6pt}] (-.24,-.34)--(.24,.34);
    }
  },
  pics/capblock/.style={
    code={
      \draw[diagram] (-.78,-.30)--(-.78,.34);
      \draw[diagram] (-.30,-.30)--(-.30,.34);
      \draw[diagram] (.10,-.30)--(.10,.02)
        .. controls (.10,.38) and (.74,.38) .. (.74,.02)
        --(.74,-.30);
      \node[blue label] at (-.78,-.58) {1};
      \node[blue label] at (-.30,-.58) {3};
      \node[blue label] at (.10,-.58) {2};
      \node[blue label] at (.74,-.58) {4};
    }
  },
  pics/nbundle/.style={
    code={
      \foreach \x in {-.27,0,.27}
        \draw[diagram] (\x,-1.7)--(\x,1.25); 
      \draw[decorate,decoration={brace,amplitude=3.5pt},line width=.9pt]
        (-.43,1.48)--(.43,1.48)
        node[midway,yshift=9pt,font=\normalsize] {$n$};
    }
  },
  pics/tk/.style={
    code={

      \draw[diagram]
        (-.25,1.58)
        .. controls (-.20,1.23) and (-.17,1.03) .. (.13,.73)
        .. controls (.38,.43) and (.38,.22) .. (.0,-.08)
        .. controls (-.38,-.38) and (-.38,-.60) .. (.0,-.91)
        .. controls (.28,-1.12) and (.35,-1.55) .. (.64,-1.55)
        .. controls (.86,-1.55) and (1.00,-1.26) .. (1.00,-.77)
        --(1.00,1.2);
        
       \draw[diagram,preaction={draw=white,line width=4.0pt}]
        (-1.00,1.58)--(-1.00,-.77)
        .. controls (-1.00,-1.26) and (-.86,-1.55) .. (-.64,-1.55)
        .. controls (-.35,-1.55) and (-.28,-1.12) .. (-.0,-.91)
        .. controls (.38,-.60) and (.38,-.38) .. (-.0,-.08)
        .. controls (-.38,.22) and (-.38,.43) .. (-.13,.73)
        .. controls (.17,1.03) and (.24,1.13) .. (.25,1.2);
        
        \draw[diagram,preaction={draw=white,line width=4.0pt}]
		(.13,.73)
		.. controls (.38,.43) and (.38,.22) .. (.0,-.08);
		\draw[diagram]
		(.06,.81) -- (.13,.73)
		.. controls (.38,.43) and (.38,.22) .. (.0,-.08)
		-- (-.1,-.165);
		
		\draw[diagram]
		(.25,1.2)
		.. controls (.3,1.6) and (.95,1.6) .. (1,1.2);

      \node[blue label] at (-1.28,1.1) {1};
      \node[blue label] at (-.43,1.1) {3};
      \node[blue label] at (.43,1.1) {2};
      \node[blue label] at (1.28,1.1) {4};
      \node[blue label] at (-1.28,-1.1) {1};
      \node[blue label] at (-.43,-1.1) {2};
      \node[blue label] at (.43,-1.1) {3};
      \node[blue label] at (1.28,-1.1) {4};
    }
  },
  pics/twoannuli/.style={
    code={
      \draw[diagram] (0,-.20)
        ellipse [x radius=1.65,y radius=.82];
      \draw[diagram] (0,-.20)
        ellipse [x radius=1.05,y radius=.52];

      \foreach \x in {-.26,0,.26}
        \draw[white,line width=8pt] (\x,-.42)--(\x,1.25);
      \foreach \x in {-.26,0,.26}{
        \draw[diagram] (\x,-.42)--(\x,1.25);
        \draw[diagram] (\x,-1.7)--(\x,-1.2);
      }
      \draw[decorate,decoration={brace,amplitude=3.5pt},line width=.9pt]
        (-.42,1.48)--(.42,1.48)
        node[midway,yshift=9pt,font=\normalsize] {$n$};
    }
  },
  pics/fourannuli1234/.style={
    code={
      \draw[diagram] (0,-.20)
        ellipse [x radius=1.90,y radius=.90];
      \draw[diagram] (0,-.20)
        ellipse [x radius=1.55,y radius=.72];
      \draw[diagram] (0,-.20)
        ellipse [x radius=1.20,y radius=.55];
      \draw[diagram] (0,-.20)
        ellipse [x radius=.85,y radius=.38];

      \foreach \x in {-.26,0,.26}
        \draw[white,line width=8pt] (\x,-.3)--(\x,1.25);
      \foreach \x in {-.26,0,.26}{
        \draw[diagram] (\x,-.3)--(\x,1.25);
        \draw[diagram] (\x,-1.7)--(\x,-1.25);
      }
      \draw[decorate,decoration={brace,amplitude=3.5pt},line width=.9pt]
        (-.42,1.48)--(.42,1.48)
        node[midway,yshift=9pt,font=\normalsize] {$n$};
      \node[arc label] at (-1.8,-.2) {1};
      \node[arc label] at (-1.4,-.2) {2};
      \node[arc label] at (-1.0,-.2) {3};
      \node[arc label] at (-.6,-.2) {4};
    }
  },
  pics/fourannuliPermuted/.style={
    code={
      \draw[diagram] (0,-.20)
        ellipse [x radius=1.90,y radius=.90];
      \draw[diagram] (0,-.20)
        ellipse [x radius=1.55,y radius=.72];
      \draw[diagram] (0,-.20)
        ellipse [x radius=1.20,y radius=.55];
      \draw[diagram] (0,-.20)
        ellipse [x radius=.85,y radius=.38];

      \foreach \x in {-.26,0,.26}
        \draw[white,line width=8pt] (\x,-.3)--(\x,1.25);
      \foreach \x in {-.26,0,.26}{
        \draw[diagram] (\x,-.3)--(\x,1.25);
        \draw[diagram] (\x,-1.7)--(\x,-1.25);
      }
      \draw[decorate,decoration={brace,amplitude=3.5pt},line width=.9pt]
        (-.42,1.48)--(.42,1.48)
        node[midway,yshift=9pt,font=\normalsize] {$n$};
      \node[arc label] at (-1.8,-.2) {1};
      \node[arc label] at (-1.4,-.2) {3};
      \node[arc label] at (-1.0,-.2) {2};
      \node[arc label] at (-.6,-.2) {4};
    }
  }
}

\begin{tikzpicture}[font=\large]
  \node at (2,4.70) {$T_U\,=$};
  \pic[scale=.80] at (3.20,4.70) {cup};

  \node at (5.35,4.70) {$T_K\,=$};
  \pic at (7.55,4.65) {tk};

  \node[anchor=east] at (12.95,4.70) {$A_n(T_U):$};
  \pic at (14.20,4.70) {nbundle};
  \draw[movie arrow] (15.35,4.70)--(17.90,4.70);
  \node at (16.15,5.42) {$S^1\,\times$};
  \pic[scale=.55] at (16.95,5.42) {cup};
  \pic at (20.10,4.65) {twoannuli};

  \node[anchor=east] at (2.5,0) {$A_n(T_K):$};
  \pic at (3.15,0) {nbundle};

  \draw[movie arrow] (4.05,0)--(6.35,0);
  \node[font=\small] at (4.5,.80) {$S^1\,\times$};
  \pic[scale=.52] at (5.55,.80) {doublecup};

  \pic at (8.65,-.02) {fourannuli1234};

  \draw[movie arrow] (10.95,0)--(13.3,0);
  \node[font=\small] at (11.35,.80) {$S^1\,\times$};
  \node[font=\small] at (11.9,.80) {$\bigl($};
  \pic[scale=.48] at (12.45,.80) {crossblock};
  \node[font=\small] at (13.1,.83) {$\bigr)^{3}$};

  \pic at (15.55,-.02) {fourannuliPermuted};

  \draw[movie arrow] (17.75,0)--(20,0);
  \node[font=\small] at (18.25,.80) {$S^1\,\times$};
  \pic[scale=.52] at (19.25,.82) {capblock};

  \pic at (22.05,-.02) {twoannuli};
\end{tikzpicture}